# ON THE EXISTENCE OF REDUCED ORDER PROPORTIONAL INTEGRAL OBSERVER FOR THE STATE ESTIMATION OF CONTINUOUS-TIME LINEAR TIME-INVARIANT SYSTEMS

Konstadinos H. Kiritsis

Hellenic Air Force Academy, Department of Aeronautical Sciences, Division of Automatic Control, Dekelia Air Base, PC 13671 Acharnes, Attikis,Tatoi, Greece
e-mail: konstantinos.kyritsis@hafa.haf.gr

**Abstract**

*In this paper the explicit necessary and sufficient conditions for the existence of reduced order proportional-integral observer for the state estimation of continuous-time linear time-invariant systems are established. A procedure is given for the calculation of observer matrices.*

## 1. INTRODUCTION

In [1] Wojciechowski added an additional term to Luenberger's full order observer [2] for the state estimation of single-input single-output continuous-time linear time-invariant systems. This term is proportional to the integral of the output estimation error. The new observer was called proportional-integral observer. In [3] implicit necessary and sufficient conditions for the existence of proportional-integral observer for the state estimation of linear multivariable time-varying systems have been established. In [4] implicit necessary and sufficient conditions for the existence of reduced order proportional-integral observer for the state estimation of linear multivariable time-varying systems have been established. In [5] necessary and sufficient conditions have been established under which the proportional-integral observer achieves Exact Loop Transfer Recovery. In [6] it was proved that the proportional-integral observer can estimate the state of continuous-time linear time-invariant systems and the state of dynamical systems with arbitrary external input which appear as unknown input, nonlinearity or unmodelled dynamics

A parametric design method for proportional-integral observers for the state estimation of continuous -time linear time-invariant systems was developed in [7].

There are many approaches for fault detection using proportional-integral observers; for more complete references, we refer the reader to [8-10] and the papers cited therein. The proportional-integral observer literature is rich; for more complete references, we refer the reader to [11] and [12] and the papers cited therein. In this paper explicit necessary and sufficient conditions for the existence of a reduced order proportional-integral observer for the state estimation of continuous-time linear time-invariant systems are established. Furthermore a procedure is given for the computation of matrices of the reduced order proportional-integral observer.

## 2. PROBLEM STATEMENT

Consider a linear time-invariant system described by the following state-space equations

$$\dot{x}(t) = Ax(t) + Bu(t) \tag{1}$$

$$y(t) = Cx(t) \tag{2}$$

where $A$, $B$ and $C$ are real matrices of dimensions ($n \times n$), ($n \times m$) and ($p \times n$) respectively, $x(t)$ is the state vector of size ($n \times 1$), $u(t)$ is the vector of inputs of size ($m \times 1$) and $y(t)$ is the vector of outputs of size ($p \times 1$). In what follows, we assume without any loss of generality that

$$rank[C] = p \tag{3}$$

Relationship (3) implies the existence of a non-singular real matrix $T$ of size ($n \times n$) such that

$$CT = [I_p, 0] \tag{4}$$

Introducing the following similarity transformation

$$x(t) = Tz(t) \tag{5}$$

and taking into account (4) and (5) the state-space equations (1) and (2) of the given system can be expressed as

$$\begin{bmatrix} \dot{z}_1(t) \\ \dot{z}_2(t) \end{bmatrix} = \begin{bmatrix} A_{11} & A_{12} \\ A_{21} & A_{22} \end{bmatrix} \begin{bmatrix} z_1(t) \\ z_2(t) \end{bmatrix} + \begin{bmatrix} G_1 \\ G_2 \end{bmatrix} u(t) \tag{6}$$

$$y(t) = [I_p, 0] \begin{bmatrix} z_1(t) \\ z_2(t) \end{bmatrix} \tag{7}$$

where

$$\begin{bmatrix} A_{11} & A_{12} \\ A_{21} & A_{22} \end{bmatrix} = T^{-1}AT, \quad T^{-1}B = \begin{bmatrix} G_1 \\ G_2 \end{bmatrix} \quad \text{and} \quad z(t) = T^{-1}x(t) = \begin{bmatrix} z_1(t) \\ z_2(t) \end{bmatrix} \tag{8}$$

where $A_{11}$, $A_{12}$, $A_{21}$ and $A_{22}$ are real matrices of dimensions ($p \times p$), ($p \times (n-p)$), (($n-p) \times p$) and (($n-p) \times (n-p)$), respectively and $z_1(t)$ and $z_2(t)$ are vectors of dimensions (p x 1) and ((n–p) x 1) )), respectively. From (7) it is obvious that $\mathbf{z}_1(t)$ denotes the states that are measurable and $\mathbf{z}_2(t)$ denotes the states that are not measurable. Using (8), equations (6) and (7) can be expressed as [13]

$$\dot{\mathbf{z}}_1(t) = \mathbf{A}_{11}\mathbf{z}_1(t) + \mathbf{A}_{12}\mathbf{z}_2(t) + \mathbf{G}_1\,\mathbf{u}(t) \tag{9}$$

$$\dot{\mathbf{z}}_2(t) = \mathbf{A}_{21}\mathbf{z}_1(t) + \mathbf{A}_{22}\mathbf{z}_2(t) + \mathbf{G}_2\,\mathbf{u}(t) \tag{10}$$

$$\mathbf{y}(t) = \mathbf{z}_1(t) \tag{11}$$

From (11) we have that

$$\dot{\mathbf{y}}(t) = \dot{\mathbf{z}}_1(t) \tag{12}$$

Substituting (11) and (12) into (9) and (10) and after some algebraic manipulations we obtain [20]

$$\dot{\mathbf{z}}_2(t) = \mathbf{A}_{22}\mathbf{z}_2(t) + \mathbf{D}\,\mathbf{u}_1(t) \tag{13}$$

$$\mathbf{y}_1(t) = \mathbf{A}_{12}\mathbf{z}_2(t) \tag{14}$$

The vectors $\mathbf{u}_1(t)$, $\mathbf{y}_1(t)$ and the matrix $\mathbf{D}$ of appropriate dimensions are given by

$$\mathbf{u}_1(t) = \begin{bmatrix} \mathbf{y}(t) \\ \mathbf{u}(t) \end{bmatrix} \tag{15}$$

$$\mathbf{y}_1(t) = \dot{\mathbf{y}}(t) - \mathbf{A}_{11}\mathbf{y}(t) - \mathbf{G}_1\,\mathbf{u}(t) \tag{16}$$

$$\mathbf{D} = [\mathbf{A}_{21},\ \mathbf{G}_2\,] \tag{17}$$

Consider also a linear time-invariant system described by the following state-space equations

$$\dot{\hat{\mathbf{z}}}_2(t) = [\mathbf{A}_{22} - \mathbf{L}\mathbf{A}_{12}]\hat{\mathbf{z}}_2(t) + \mathbf{L}\mathbf{y}_1(t) + \mathbf{D}\mathbf{u}_1(t) + \mathbf{F}\boldsymbol{\omega}(t) \tag{18}$$

$$\dot{\boldsymbol{\omega}}(t) = \mathbf{G}[\mathbf{y}_1(t) - \mathbf{A}_{12}\,\hat{\mathbf{z}}_2(t)] \tag{19}$$

where $\hat{\mathbf{z}}_2(t)$ is the state vector of dimensions ((n–p) x 1), $\boldsymbol{\omega}(t)$ is a vector of dimensions (k x 1) and $\mathbf{L}$, $\mathbf{F}$ and $\mathbf{G}$ are real matrices of size ((n–p) x p), ((n–p) x k) and (k x p) respectively. The linear time-invariant system described by the equations (18) and (19) is a proportional-integral observer of order (n–p) for the system described by the equations (13) and (14) or alternatively a reduced order proportional-integral observer of order (n–p) for the system described by the equations (1) and (2), if and only if for arbitrary initial conditions $\hat{\mathbf{z}}_2(0)$, $\mathbf{z}_2(0)$ and any input $\mathbf{u}_1(t)$ the following relationships hold [3] (see also [7])

$$lim_{t \to +\infty}\ \mathbf{e}(t) = \mathbf{0} \tag{20}$$

$$lim_{t \to +\infty}\ \boldsymbol{\omega}(t) = \mathbf{0} \tag{21}$$

where $\mathbf{e}(t) = [\hat{\mathbf{z}}_2(t) - \mathbf{z}_2(t)]$ is the state estimation error, $\hat{\mathbf{z}}_2(t)$ is an estimate of the state vector $\mathbf{z}_2(t)$ and $\boldsymbol{\omega}(t)$ is a vector representing the integral of the weighted output estimation error [7]. The relationships (20) and (21) are simultaneously satisfied if and only if the following matrix

$$\begin{bmatrix} A_{22} - LA_{12} & F \\ -GA_{12} & 0 \end{bmatrix} \quad (22)$$

of dimensions (($n-p+k$) x ($n-p+k$)) is Hurwitz stable, i.e. all its eigenvalues have negative real parts [3] (see also [7]). Thus the problem of the design of the reduced order proportional-integral observer of order ($n-p$) can be stated as follows: Do there exist real matrices **L**, **F** and **G** of appropriate dimensions such that the matrix given by (22) is Hurwitz stable? If so, give conditions for existence and a procedure for the calculation of the real matrices **L**, **F** and **G**.

## 2. BASIC CONCEPTS AND PRELIMINARY RESULTS

Let $C$ be the field of complex numbers, also let $C^+$ be the set of all complex numbers $\lambda$ with $Re(\lambda) \geq 0$.

Let **A** be a nonzero matrix over $C$ of size ($p$ x $q$). The row rank (column rank) of a matrix **A** is the maximum value $r$ for which there exist $r$ linearly independent rows (columns) of **A**. The rank, the row rank and the column rank of **A** are all equal [14]. The rank of the matrix **A** cannot exceed the minimum of $p$ and $q$, that is

$$rank[\mathbf{A}] \leq min\ (p, q) \quad (23)$$

Let **A** and **B** be matrices of size ($p$ x $m$) and ($m$ x $q$) respectively. Then [15]

$$rank[\mathbf{AB}] \leq rank[\mathbf{B}]; \quad rank[\mathbf{AB}] \leq rank[\mathbf{A}] \quad (24)$$

**Definition 1:** The polynomial $c(s)$ with real coefficients is said to be strictly Hurwitz if and only if $c(s) \neq 0, \forall s \in C^+$.

**Definition 2:** The real matrix **A** of dimensions ($q$ x $q$) is said to be Hurwitz stable if and only if and only if the characteristic polynomial of the matrix **A** is a strictly Hurwitz polynomial

Let **A** be a real matrix of size ($q$ x $q$). The spectrum of the real matrix **A**, is the set of all its eigenvalues and is denoted by $\sigma(\mathbf{A})$.

**Lemma 1:** Let **A** and **C** be real matrices of dimensions ($n$ x $n$) and ($p$ x $n$), respectively and **C** not zero. Further let $\sigma(\mathbf{A})$ be the spectrum of the matrix **A**. Then the pair (**A, C**) is detectable [16] if and only if one of the following equivalent conditions holds [17]:

(a) $rank \begin{bmatrix} \mathbf{C} \\ \mathbf{I}_n s - \mathbf{A} \end{bmatrix} = n$, $\forall s \in C^+$

(b) $rank \begin{bmatrix} \mathbf{C} \\ \mathbf{I}_n \lambda - \mathbf{A} \end{bmatrix} = n$, $\forall \lambda \in \sigma(\mathbf{A})$ with $Re(\lambda) \geq 0$

**Lemma 2.** Let $\mathbf{A}$ and $\mathbf{C}$ be real matrices of dimensions *(n x n)* and *(p x n)*, respectively. Then the pair $(\mathbf{A}, \mathbf{C})$ is observable if and only if for every monic polynomial $c(s)$ over the ring of polynomials with real coefficients of degree $n$, there exists a real matrix $\mathbf{K}$ of dimensions *( n x p)*, such that the matrix $[\mathbf{A}+\mathbf{KC}]$ has characteristic polynomial $c(s)$[18].

**Lemma 3:** Let $\mathbf{A}$ be a real matrix of dimensions *(n x n)*. Then the matrix $\mathbf{A}$ is Hurwitz stable if and only if the following condition holds:

(a) $rank[\mathbf{I}_n s - \mathbf{A}] = n$, $\forall s \in C^+$

*Proof:* Let $\mathbf{A}$ be a Hurwitz stable real matrix of dimensions *(n x n)*. From Definition 2 it follows that the characteristic polynomial $c(s)$ of the real matrix $\mathbf{A}$ is a strictly Hurwitz polynomial and therefore from Definition 1 it follows that

$$c(s) \neq 0, \forall s \in C^+ \quad (25)$$

The Smith-McMillan form of polynomial matrix $[\mathbf{I}_n s - \mathbf{A}]$ over the ring of polynomials with real coefficients is given by [18]

$$\mathbf{K}(s)[\mathbf{I}_n s - \mathbf{A}]\mathbf{L}(s) = [diag[c_1(s), c_2(s), ...., c_n(s)]] \quad (26)$$

where $\mathbf{K}(s)$ and $\mathbf{L}(s)$ are unimodular matrices over the ring of polynomials with real coefficients [18]. The polynomials $c_i(s)$ for $i=1, 2,...,n$ are the invariant polynomials of the real matrix $\mathbf{A}$ and therefore [18]

$$c(s) = \Pi_{i=1}^{n} c_i(s) \quad (27)$$

From (25) and (27) we have that

$$c_i(s) \neq 0, \forall s \in C^+, \forall i=1,2,...,n \quad (28)$$

From (28) we obtain:

$$rank\{diag[c_1(s), c_2(s), ...., c_n(s)]\} = n, \forall s \in C^+ \quad (29)$$

Since $\mathbf{K}(s)$ and $\mathbf{L}(s)$ are unimodular matrices over the ring of polynomials with real coefficients [18], from (26) we obtain [18]:

$$rank[\mathbf{I}_n s - \mathbf{A}] = rank\{diag[c_1(s), c_2(s), ...., c_n(s)]\} \quad (30)$$

Relationships (29) and (30) imply that

$$rank[\mathbf{I}_n s - \mathbf{A}] = n, \forall s \in C^+ \quad (31)$$

This is condition (a) of Lemma 3. To prove sufficiency, we assume that condition (a) of Lemma 3 holds. Since by assumption condition (a) holds we have that

$$rank[\mathbf{I}_n s - \mathbf{A}] = n \ , \forall s \in C^+ \tag{32}$$

From (32) and (30) we obtain:

$$rank\{diag[c_1(s), c_2(s), ...., c_n(s)]\} = n, \forall s \in C^+ \tag{33}$$

From (33) we have that

$$c_i(s) \neq 0 \ , \forall s \in C^+ \ , \ \forall i=1,2,...,n \tag{34}$$

From (34) we have that

$$\prod_{i=1}^{n} c_i(s) \neq 0 \ , \forall s \in C^+, \ \forall i=1,2,...,n \tag{35}$$

Relationships (35) and (27) imply that

$$c(s) \neq 0 \ , \forall s \in C^+ \tag{36}$$

Relationship (36) and Definition 1 imply that $c(s)$ is a strictly Hurwitz polynomial. Since by assumption $c(s)$ is the characteristic polynomial of the real matrix $\mathbf{A}$, from Definition 2 it follows that the real matrix $\mathbf{A}$ is Hurwitz stable. This completes the proof.

The following Lemma is taken from [19] (see also [20]).

**Lemma 4:** Let $\mathbf{A}$ and $\mathbf{C}$ be real matrices of dimensions *(n x n)* and *(p x n)*, respectively and $\mathbf{C}$ not zero. Further, let the pair $(\mathbf{A}, \mathbf{C})$ be detectable. Then there exists a real matrix $\mathbf{K}$ of dimensions *(n x p)* such that the real matrix $[\mathbf{A}+\mathbf{KC}]$ is Hurwitz stable.

The following Lemma is taken from [19] (see also [20]).

**Lemma 5:** Let $\mathbf{A}$ and $\mathbf{C}$ be real matrices of size *(n x n)*, *(p x n)* respectively as in (1) and (2). Further, let $\mathbf{A}_{12}$ and $\mathbf{A}_{22}$ be real matrices of appropriate dimensions as in (8). Then the pair $(\mathbf{A}, \mathbf{C})$ is detectable if and only if the pair $(\mathbf{A}_{22}, \mathbf{A}_{12})$ is detectable.

## 4. MAIN RESULTS

The Theorems that follow are the main results of this paper.

**Theorem 1.** Let $rank[\mathbf{A}_{12}] = q$. Then the system described by equations (18) and (19) is a reduced order proportional-integral observer for the system described by equations (1) and (2), only if the following condition holds:

(a) $rank[\mathbf{G}] = k$

(b) $q \geq k$

***Proof***: Let the system described by equations (18) and (19) is a reduced order proportional-integral observer for the system described by equations (1) and (2). Then the real matrix of dimensions $((n-p+k) \times (n-p+k))$ given by (22) is Hurwitz stable. Hurwitz stability of the matrix given by (22) and Lemma 3 imply that

$$rank\begin{bmatrix} \mathbf{I}_{n-p}s - \mathbf{A}_{22} + \mathbf{LA}_{12} & -\mathbf{F} \\ \mathbf{GA}_{12} & \mathbf{I}_k s \end{bmatrix} = (n-p+k) \ , \ \forall s \in C^+ \tag{37}$$

Since the $(n - p + k)$ rows of the matrix on the left side of (37) are linearly independent over $C$, $\forall s \in C^+$, a subset of these columns consisting of the last $k$ rows must be also linearly independent over $C$, $\forall s \in C^+$; therefore

$$rank[\mathbf{GA}_{12}, \mathbf{I}_k s] = k \ , \ \forall s \in C^+ \tag{38}$$

For $s=0$ from (38) we obtain:

$$rank[\mathbf{GA}_{12}, \mathbf{0}] = rank[\mathbf{GA}_{12}] = k \tag{39}$$

From (39) and (24) it follows that

$$rank[\mathbf{GA}_{12}] = k \leq rank[\mathbf{G}] \tag{40}$$

Since the matrix **G** is of size $(k \times p)$, from (40) and (23) we obtain:

$$rank[\mathbf{G}] = k \tag{41}$$

This is condition (a) of Theorem 1. From (40) and (24) it follows that

$$rank[\mathbf{GA}_{12}] = k \leq rank[\mathbf{A}_{12}] \tag{42}$$

Since by assumption $rank[\mathbf{A}_{12}] = q$, condition (b) of Theorem 1 follows from (42) and the proof is complete.

**Theorem 2.** Let $rank[\mathbf{A}_{12}] = q$. Then the system described by equations (18) and (19) with $rank[\mathbf{G}] = k$ and $q \geq k$ is a reduced order proportional-integral observer for the system described by equations (1) and (2), if and only if the following condition holds:

(a) The pair (**A**, **C**) is detectable.

***Proof***: Let the system described by equations (18) and (19) with $rank[\mathbf{G}] = k$ and $q \geq k$ is a reduced proportional-integral observer for the system described by equations (1) and (2). Then the real matrix of dimensions $((n-p+k) \times (n-p+k))$ given by (22) is Hurwitz stable. Hurwitz stability of the matrix given by (22) and Lemma 3 imply that

$$rank\begin{bmatrix} \mathbf{I}_{n-p}s - \mathbf{A}_{22} + \mathbf{LA}_{12} & -\mathbf{F} \\ \mathbf{GA}_{12} & \mathbf{I}_k s \end{bmatrix} = (n-p+k) \ , \ \forall s \in C^+ \tag{43}$$

Since the $(n - p + k)$ columns of the matrix on the left side of (43) are linearly independent over $C$, $\forall s \in C^+$, a subset of these columns consisting of the first $(n-p)$ columns must be also linearly independent over $C$, $\forall s \in C^+$; therefore

$$rank\begin{bmatrix} \mathbf{I}_{n-p}s - \mathbf{A}_{22} + \mathbf{LA}_{12} \\ \mathbf{GA}_{12} \end{bmatrix} = rank\begin{bmatrix} \mathbf{GA}_{12} \\ \mathbf{I}_{n-p}s - \mathbf{A}_{22} + \mathbf{LA}_{12} \end{bmatrix} = (n - p), \forall s \in C^+ \quad (44)$$

From (44) and after simple algebraic manipulations we obtain:

$$rank\begin{bmatrix} \mathbf{GA}_{12} \\ \mathbf{I}_{n-p}s - \mathbf{A} + \mathbf{LA}_{12} \end{bmatrix} = rank\{\begin{bmatrix} \mathbf{G} & \mathbf{0} \\ \mathbf{L} & \mathbf{I}_{n-p} \end{bmatrix} \begin{bmatrix} \mathbf{A}_{12} \\ \mathbf{I}_{n-p}s - \mathbf{A}_{22} \end{bmatrix}\} = (n - p), \forall s \in C^+ \quad (45)$$

Since the matrix

$$\begin{bmatrix} \mathbf{A}_{12} \\ \mathbf{I}_{n-p}s - \mathbf{A}_{22} \end{bmatrix} \quad (46)$$

is a complex matrix for any $s \in C^+$, from (46) and (24) it follows that

$$(n - p) = rank\{\begin{bmatrix} \mathbf{G} & \mathbf{0} \\ \mathbf{L} & \mathbf{I}_{n-p} \end{bmatrix} \begin{bmatrix} \mathbf{A}_{12} \\ \mathbf{I}_{n-p}s - \mathbf{A}_{22} \end{bmatrix}\} \leq rank\begin{bmatrix} \mathbf{A}_{12} \\ \mathbf{I}_{n-p}s - \mathbf{A}_{22} \end{bmatrix}, \forall s \in C^+ \quad (47)$$

Since the complex matrix given by (46) is of size $(n \times (n-p))$, from (47) and (23) we obtain:

$$rank\begin{bmatrix} \mathbf{A}_{12} \\ \mathbf{I}_{n-p}s - \mathbf{A}_{22} \end{bmatrix} = (n - p), \forall s \in C^+ \quad (48)$$

Relationship (48) and condition (a) of Lemma 1 imply detectability of the pair $(\mathbf{A}_{22}, \mathbf{A}_{12})$. Detectability of the pair $(\mathbf{A}_{22}, \mathbf{A}_{12})$ and Lemma 5 imply detectability of the pair $(\mathbf{A}, \mathbf{C})$. This is condition (a) of Theorem 2.

To prove sufficiency, we assume that condition (a) of Theorem 2 holds. Detectability of the pair $(\mathbf{A}, \mathbf{C})$ and Lemma 5 imply detectability of the pair $(\mathbf{A}_{22}, \mathbf{A}_{12})$. Detectability of the pair $(\mathbf{A}_{22}, \mathbf{A}_{12})$ and Lemma 4 imply the existence of a real matrix $\mathbf{K}$ of dimensions $((n-p) \times p)$ such that the matrix $[\mathbf{A}_{22} + \mathbf{KA}_{12}]$ is Hurwitz stable, that is

$$det[\mathbf{I}_{n-p}s - \mathbf{A}_{22} - \mathbf{KA}_{12}] = c(s) \quad (49)$$

where $c(s)$ is a strictly Hurwitz polynomial with real coefficients of degree $(n - p)$. The real matrix $\mathbf{K}$ in (49) can be calculated using known methods for the solution of pole assignment problem by state feedback [18] (see also [19-20]). $rank[\mathbf{A}_{12}] = q$ implies the existence of non-singular real matrices $\mathbf{P}$ and $\mathbf{Q}$ of appropriate dimensions such that

$$\mathbf{A}_{12} = \mathbf{P}\begin{bmatrix} \mathbf{I}_q & \mathbf{0} \\ \mathbf{0} & \mathbf{0} \end{bmatrix}\mathbf{Q} \quad (50)$$

We set

$$\mathbf{G}=[\mathbf{I}_k,\ \mathbf{0}]\ \mathbf{P}^{-1}, \qquad q \geq k \tag{51}$$

Let $\mathbf{\Phi}$ be an arbitrary Hurwitz stable real matrix of dimensions ($k \times k$). Furthermore, let $\mathbf{X}$ be a real matrix of dimensions (($n-p$) x $k$) given by

$$\mathbf{X}=\mathbf{Q}^{-1}\begin{bmatrix}-\mathbf{\Phi}\\ \mathbf{0}\end{bmatrix} \tag{52}$$

From (50), (51) and (52) we obtain:

$$-\mathbf{GA}_{12}\mathbf{X}= [-\mathbf{I}_k, \mathbf{0}]\mathbf{P}^{-1}\mathbf{P}\begin{bmatrix}\mathbf{I}_q & \mathbf{0}\\ \mathbf{0} & \mathbf{0}\end{bmatrix}\mathbf{Q}\ \mathbf{Q}^{-1}\begin{bmatrix}-\mathbf{\Phi}\\ \mathbf{0}\end{bmatrix}=\mathbf{\Phi} \tag{53}$$

Now we form the real matrix $\mathbf{M}$ of dimensions (($n-p+k$) x ($n-p+k$)) [21]

$$\mathbf{M} = \begin{bmatrix}\mathbf{I}_{n-p} & \mathbf{X}\\ \mathbf{0} & \mathbf{I}_k\end{bmatrix}$$

The real matrix $\mathbf{M}$ is non-singular and its inverse is given by

$$\mathbf{M}^{-1}=\begin{bmatrix}\mathbf{I}_{n-p} & -\mathbf{X}\\ \mathbf{0} & \mathbf{I}_k\end{bmatrix}$$

We have that

$$\mathbf{M}^{-1}\begin{bmatrix}\mathbf{A}_{22}-\mathbf{LA}_{12} & \mathbf{F}\\ -\mathbf{GA}_{12} & \mathbf{0}\end{bmatrix}\mathbf{M} =$$

$$= \begin{bmatrix}\mathbf{A}_{22}+(-\mathbf{L}+\mathbf{XG})\mathbf{A}_{12} & (\mathbf{A}_{22}-\mathbf{LA}_{12})\mathbf{X}+\mathbf{F}+\mathbf{XGA}_{12}\mathbf{X}\\ -\mathbf{GA}_{12} & -\mathbf{GA}_{12}\mathbf{X}\end{bmatrix} \tag{54}$$

Furthermore, we set:

$$\mathbf{L}=\mathbf{XG}-\mathbf{K} \tag{55}$$

$$\mathbf{F}=-(\mathbf{A}_{22}-\mathbf{LA}_{12})\mathbf{X}-\mathbf{XGA}_{12}\mathbf{X} \tag{56}$$

Now by substituting (55), (56) and (53) into (54) we obtain:

$$\mathbf{M}^{-1}\begin{bmatrix}\mathbf{A}_{22}-\mathbf{LA}_{12} & \mathbf{F}\\ -\mathbf{GA}_{12} & \mathbf{0}\end{bmatrix}\mathbf{M} = \begin{bmatrix}\mathbf{A}_{22}+\mathbf{KA}_{12} & \mathbf{0}\\ -\mathbf{GA}_{12} & \mathbf{\Phi}\end{bmatrix} \tag{57}$$

Since by (49) the matrix $[\mathbf{A}_{22}+\mathbf{KA}_{12}]$ is Hurwitz stable and the real matrix $\mathbf{\Phi}$ is by assumption Hurwitz stable, from Lemma 3 it follows that

$$rank[\mathbf{I}_{n-p}s-\mathbf{A}_{22}-\mathbf{KA}_{12}] = (n-p), \forall s \in C^+ \tag{58}$$

$$rank[\mathbf{I}_k s-\mathbf{\Phi}] = k, \forall s \in C^+ \tag{59}$$

From (58) and (59) we obtain:

$$rank\begin{bmatrix}\mathbf{I}_{n-p}s-\mathbf{A}_{22}-\mathbf{KA}_{12} & \mathbf{0}\\ \mathbf{GA}_{12} & \mathbf{I}_k s-\mathbf{\Phi}\end{bmatrix} = (n-p+k), \forall s \in C^+ \tag{60}$$

Since by (57) the matrices

$$\begin{bmatrix}\mathbf{A}_{22}-\mathbf{LA}_{12} & \mathbf{F}\\ -\mathbf{GA}_{12} & \mathbf{0}\end{bmatrix}, \begin{bmatrix}\mathbf{A}_{22}+\mathbf{KA}_{12} & \mathbf{0}\\ -\mathbf{GA}_{12} & \mathbf{\Phi}\end{bmatrix} \tag{61}$$

of size $((n-p+k) \ x \ (n-p+k))$ are similar, we obtain:

$$rank\begin{bmatrix} \mathbf{I}_{n-p}s - \mathbf{A}_{22} + \mathbf{LA}_{12} & -\mathbf{F} \\ \mathbf{GA}_{12} & \mathbf{I}_k s \end{bmatrix} = rank\begin{bmatrix} \mathbf{I}_{n-p}s - \mathbf{A}_{22} - \mathbf{KA}_{12} & 0 \\ \mathbf{GA}_{12} & \mathbf{I}_k s - \mathbf{\Phi} \end{bmatrix} \quad (62)$$

From (60) and (62), we obtain:

$$rank\begin{bmatrix} \mathbf{I}_{n-p}s - \mathbf{A}_{22} + \mathbf{LA}_{12} & -\mathbf{F} \\ \mathbf{GA}_{12} & \mathbf{I}_k s \end{bmatrix} = (n - p + k), \forall s \in C^+ \quad (63)$$

Furthermore, from (63) and Lemma 3 it follows that matrix

$$\begin{bmatrix} \mathbf{A}_{22} - \mathbf{LA}_{12} & \mathbf{F} \\ -\mathbf{GA}_{12} & 0 \end{bmatrix}$$

with **L**, **F** and **G** given by (55), (56) and (51) respectively is Hurwitz stable and therefore according to (22) the system described by equations (18) and (19) with **L**, **F** and **G** given by (55), (56) and (51) respectively is a reduced order proportional-integral observer for the system described by equations (1) and (2). This completes the proof.

The proof of Theorem 2 is constructive and furnishes a procedure for the computation of matrices **L**, **F** and **G** of reduced order proportional-integral observer for the state estimation of the system described by equations (1) and (2). The major steps of this procedure are summarized below.

**CONSTRUCTION**

*Given*: **A, B, C**, $\mathbf{A}_{12}$, $\mathbf{A}_{22}$, $rank[\mathbf{A}_{12}] = q$ and $q \geq k$

*Find*: **L**, **F** and **G** with $rank[\mathbf{G}] = k$

*Step 1:* Check condition (a) of Theorem 2. If this condition is satisfied go to *Step 2*. If condition (a) of Theorem 2 is not satisfied, the construction of reduced order proportional-integral observer is impossible.

*Step 2:* Delectability of the pair (**A, C**) and Lemma 5 imply detectability of the pair ($\mathbf{A}_{22}$, $\mathbf{A}_{12}$). Detectability of the pair ($\mathbf{A}_{22}$, $\mathbf{A}_{12}$) and Lemma 4 imply the existence of real matrix **K** of dimensions $((n-p) \ x \ p)$ such that the matrix $[\mathbf{A}_{22} + \mathbf{KA}_{12}]$ is Hurwitz stable. The matrix **K** can be calculated using known methods for the solution of pole assignment problem by state feedback [18] (see also [19-20]).

*Step 3:* Find non-singular real matrices **P** and **Q** of appropriate dimensions such that

$$\mathbf{A}_{12} = \mathbf{P}\begin{bmatrix} \mathbf{I}_q & 0 \\ 0 & 0 \end{bmatrix}\mathbf{Q}$$

*Step 4:* Let **Φ** be an arbitrary Hurwitz stable real matrix of dimensions $(k \ x \ k)$. Put

$$X = Q^{-1} \begin{bmatrix} -\Phi \\ 0 \end{bmatrix}$$

$$L = XG - K$$

$$F = -(A_{22} - LA_{12})X - XGA_{12}X$$

$$G = [I_k, \ 0] \ P^{-1}$$

## 5. CONCLUSIONS

In this paper, the problem of existence of reduced order proportional-integral observer for the state estimation of continuous-time linear time-invariant systems is studied. In particular, explicit necessary and sufficient conditions for the existence of a reduced order proportional-integral observer for the state estimation of continuous-time linear time-invariant systems are established. Furthermore a procedure is given for the computation of matrices of reduced order proportional-integral observer.